\documentclass{amsart}
\usepackage{amsmath,amssymb,amsthm,amscd, mathtools, latexsym}
\usepackage{enumerate,varioref, dsfont, fancyhdr}
\usepackage{tikz-cd}
\usepackage{enumitem}
\usepackage{multicol}
\usepackage{hyperref}
\usepackage{url}
\usepackage{textcomp}

\newtheorem{Thm}{Theorem}[section]
\newtheorem{Lem}[Thm]{Lemma}
\newtheorem*{Lem*}{Lemma}
\newtheorem{Prop}[Thm]{Proposition}

\newtheorem{Cor}[Thm]{Corollary}

\newtheorem{Rem}[Thm]{Remark}

\newtheorem{Set}[Thm]{Set-up}
\newtheorem{Def}[Thm]{Definition}
\newtheorem{claim}[Thm]{Claim}

\newtheorem{Ass}[Thm]{Assumption}
{\theoremstyle{plain}

\newtheorem{Not}[Thm]{Notation}

}

\newcommand\restr[2]{\ensuremath{\left.#1\right|_{#2}}}

\def\rit{{\mathbb R}}
\def\cit{{\mathbb C}}

\def\qit{{\mathbb Q}}
\def\zit{{\mathbb Z}}

\def\pit{{\mathbb P}}

\def\0{{\mathcal O}}

\def\Aut{\mathop{\rm Aut}\nolimits}
\def\End{\mathop{\rm End}\nolimits}
\def\A{{\mathfrak A}}

\def\C{{\mathcal C}}

\def\G{{\mathcal G}}
\def\M{{\mathcal M}}

\def\A{{\mathcal A}}

\def\X{{\mathcal X}}
\def\Y{{\mathcal Y}}

\title{Non-divisible cycles on products of very general Abelian varieties}
\author{Humberto A. Diaz}
\newcommand{\Addresses}{{\bigskip \footnotesize
\textsc{Department of Mathematics, Washington University, St. Louis, MO 63130} \par \nopagebreak
\textit{Email address}: \ \texttt{humberto@wustl.edu}}}
\date{}
\begin{document}
\begin{abstract}
In this paper, we give a recipe for producing infinitely many non-divisible codimension $2$ cycles on a product of two or more very general Abelian varieties. In the process, we introduce the notion of ``field of definition" for cycles in the Chow group modulo (a power of) a prime. We show that for a quite general class of codimension $2$ cycles we call ``primitive cycles," the field of definition is a ramified extension of the function field of a modular variety. This ramification allows us to use Nori's isogeny method (modified by Totaro) to produce infinitely many non-divisible cycles. As an application, we prove the Chow group modulo a prime of a product of $3$ or more very general elliptic curves is infinite, generalizing work of Schoen.
\end{abstract}
\maketitle

\section{Introduction}
\noindent A fairly recent result of Totaro \cite{T} shows that there exist smooth complex projective varieties whose Chow group modulo $\ell$ is infinite for all primes $\ell$, generalizing earlier work of Rosenschon and Srinivas in \cite{RS}. Equivalently, there are examples of such varieties for which the Chow group is far from being an $\ell$-divisible group for any prime $\ell$. The literature gives a general recipe for proving this type of result. Indeed, the methods employed in \cite{RS} and \cite{T} show the non-divisibility of the Ceresa cycle on the very general Abelian threefold. Then, an adaptation of Nori's isogeny technique is used to produce the required infinitely many cycles. (Note that in the case of \cite{RS}, the authors use Atkin-Lehner correspondences as in \cite{S1}). \\
\indent Our goal here will be to extend the existing techniques in the literature for producing infinitely many cycles modulo $\ell$ to prove similar results for products of two or more very general Abelian varieties. In particular, we are able to use the method developed here to prove the following result, which extends work of Schoen in \cite{S1}:
\begin{Thm}\label{main} Let $E_{1}$, $E_{2}, \ldots E_{n}$ be very general complex elliptic curves. Then,
$$CH^{d} (E_{1}\times E_{2} \times \dots \times E_{n}) \otimes \zit/\ell$$
is infinite for all primes $\ell$, $2 \leq d \leq n-1$ and $n\geq 3$.
\end{Thm}
\noindent As noted in the sequel, it will be sufficient to show this for $d=2$ and $n=3$. We briefly describe the method involved, which works best in codimension $2$ (because of the Merkurjev-Suslin theorem). In general, we consider a family of smooth complex projective varieties $f: \X \to S$ whose fibers are products of two or more very general Abelian varieties. We introduce the notion of a {\em primitive correspondence} $P$, which is a correspondence on $\X$ which kills the first level of the coniveau filtration on $H^{3}$ of the very general fiber of $f$. A codimension $2$ cycle on $\X$ is said to be {\em primitive with respect to $P$} if the normal function of $P_{*}(\gamma)$ is non-torsion. Our first order of business will be to use the techniques of \cite{RS} and \cite{T} to show that any primitive cycle $\gamma$ is non-divisible in the Chow group of the very general fiber of $f$. In the case of Theorem \ref{main}, $\gamma$ is an ersatz Ceresa cycle which was considered in \cite{BST} and \cite{S1}. We spread out this cycle to a larger family $f: \X \to S$ and using some rather elementary geometric observations about the total space $\X$, we are able to deduce the required statement about the normal function of $\gamma$. (We observe that in the case of \cite{RS} and \cite{T}, the total space of the universal Abelian threefold with level structure is much more complicated and some deep results of Hain \cite{H} must be invoked.)\\
\indent In order to produce infinitely many cycles, we use Nori's isogeny method. To implement his method, we introduce the notion of a ``field of definition" for a cycle in the Chow group modulo a power of $\ell$. This is defined as the minimum field for which the absolute Galois group acts trivially on $\gamma$. (In general, this notion does not seem to be well-understood and reflects the general failure of Galois descent for cycles modulo $\ell$). We will need to show that the field of definition of $\gamma$ is a ramified extension of a product of Siegel modular varieties (the infinitude is then detected by the action of isogenies on the ramification locus). This should not be surprising; indeed, as we show in the sequel, primitive cycles cannot be defined on (products of) universal Abelian varieties. We note that in the case of the Ceresa cycle on the very general Abelian threefold \cite{N}, \cite{RS} and \cite{T}, the field of definition and the ramification locus can be described rather explicitly (see Remark \ref{definition}). In general, however, it is hard to give such an explicit characterization.
\begin{Rem} For a smooth complex projective variety, the Chow group of cycles modulo algebraic equivalence is a divisible group. Hence, Theorem \ref{main} gives yet another instance of a variety whose Griffiths group is of infinite rank. The first example of this phenomenon was provided by Clemens \cite{C}, who used techniques different from the isogeny technique of Nori.
\end{Rem}
\indent As a consequence of Theorem \ref{main}, one obtains the result below (as in \cite{RS} Theorem 1.4 and \cite{T} Corollary 0.2) using Schoen's external product maps (\cite{S3} Theorem 0.2), which also involves very general elliptic curves:
\begin{Cor} Let $E_{1}$, $E_{2}, \ldots E_{n}$ be very general complex elliptic curves. Then, the $\ell$-torsion
\[CH^{d} (E_{1}\times E_{2} \times \dots \times E_{n})[\ell]\]
is infinite for all primes $\ell$, $3 \leq d \leq n-1$ and $n\geq 4$.
\end{Cor}

\subsection*{Convention} While the term ``very general" is undoubtedly understood by many readers, we give a brief definition. In this paper, whenever $f: \X \to S$ is a family of varieties over $\cit$, a ``very general" member refers to a fiber $\X_{s} = f^{-1}(s)$ for $s \in S(\cit)$ outside of a countable union of Zariski closed subsets of $S$. Since all very general fibers of such a family are non-canonically isomorphic as varieties over $\cit$, we will occasionally use the term ``the very general" member.

\subsection*{Acknowledgements}

The author would like to thank Burt Totaro for his interest in early drafts of this paper. He would also like to thank Chad Schoen for some useful comments and Robert Laterveer for calling his attention to a certain method of Claire Voisin \cite{V} described in the sequel. The author would also like to thank the referee for making many useful suggestions.

\section{The Bloch-Esnault criterion}
\noindent We begin with the following deep result which is central to the arguments that follow and which is a fixture in arguments about non-divisibility of Chow groups (c.f., \cite{S1},\cite{RS}, \cite{T}). 
\begin{Thm}[Bloch-Esnault, \cite{BE} Theorem 1.2] Let $X$ be a smooth projective variety over a complete discrete valuation field $K$ with perfect residue field $k$ of mixed characteristic $(0,p)$. Suppose further that $X$ has good ordinary reduction and that 
\begin{enumerate}[label=(\alph*)]
\item\label{no-tors-2} The crystalline cohomology of the special fiber $Y$ has no torsion.
\item\label{forms} $H^{0} (Y, \Omega^{m}_{Y}) \neq 0$.
\end{enumerate}
Then, $N^{1}H^{m}_{\text{\'et}} (X_{\overline{K}}, \zit/p) \neq H^{m}_{\text{\'et}} (X_{\overline{K}}, \zit/p)$, where $N^{*}$ denotes the coniveau filtration in \'etale cohomology.
If the \'etale cohomology $H^{m+1}_{\text{\'et}} (X_{\overline{K}}, \zit_{p})$ has no torsion, then for all $r\geq 1$,
\[ N^{1}H^{m}_{\text{\'et}} (X_{\overline{K}}, \zit/p^{r}) \neq H^{m}_{\text{\'et}} (X_{\overline{K}}, \zit/p^{r})\]
\end{Thm}
\noindent Note that the original result only gives the first statement, but the second statement may easily be deduced from a coefficients theorem. There is no torsion in the crys-
talline cohomology of an Abelian variety, so we have the following consequence:
\begin{Cor} Let $A_{1}, A_{2}, \ldots A_{s}$ be very general complex Abelian varieties of dimension $g_{1}, \ldots g_{s}$; set $A= A_{1} \times \ldots \times A_{s}$ and $g:= g_{1}+\ldots +g_{s}$. Then, 
\[ N^{1}H^{m}_{\text{\'et}} (A, \zit/p^{r}) \neq H^{m}_{\text{\'et}} (A, \zit/p^{r})\]
for all $1\leq m \leq g$, primes $p$ and $r\geq 1$.
\begin{proof} Fix a prime $p.$ Given the Bloch-Esnault criterion, it suffices to establish that the very general complex Abelian variety of dimension $g$ admits a model over a complete discrete valuation field $K$ with perfect residue field $k$ of mixed characteristic $(0,p)$ that has good ordinary reduction. A reference for this   fact is given by \cite{NO} Theorem 3.1.
\end{proof}
\end{Cor}
\section{Primitive correspondences}
\begin{Not}\label{notation} In the sequel, we will denote the absolute Galois group of a field $F$ by $G_{F}$. Moreover, we will consider smooth varieties $X$ over $\cit$ and for any Abelian group $G$ we let $H^{*} (X, G)$ be the singular cohomology group with coefficients in $G$ of the underlying complex manifold. Moreover, let $N^{*}H^{*}(X, G)$ denote the corresponding coniveau filtration. For $S$ a smooth quasi-projective variety over $\cit$ and $f: \mathcal{X} \to S$ a smooth projective morphism, let $f^{-1} (s) = \mathcal{X}_{s}$ for $s \in S(\cit)$. Moreover, let $\cit(S)$ be the function field of $S$, $\X_{\cit(S)}$ the generic fiber and for any field extension $L/\cit(S)$, let $\X_{L}:= \X_{\cit(S)} \times_{\cit(S)} L$.\\
\indent We let $\M_{S}$ denote the category of relative Chow motives over $S$ with integral coefficients, and for any smooth projective $S$-variety $\X/S$, denote its corresponding Chow motive by $\M(\X)$. Moreover, for $i\geq 0$ let $R^{i}f_{*}G$ denote the local system whose stalk over $s \in S(\cit)$ is $H^{i} (\mathcal{X}_{s}, G)$.
\end{Not}
\noindent With the above notation, consider a relative correspondence 
\[ P \in CH^{g} (\mathcal{X} \times_{S} \mathcal{X}) = \End_{\M_{S}} (M(\mathcal{X})) \]
\begin{Def} We say that $P$ is a {\em primitive correspondence in codimension $i$} if the following hold:
\begin{enumerate}[label=(\alph*)]
\item For some integer $c \geq 1$, $P^{\circ 2} = c\cdot P$.
\item\label{coniveau} For the very general $s \in S(\cit)$, the complementary correspondence $Q = c\cdot\Delta_{\mathcal{X}} - P$ satisfies
\[ Q_{*}H^{2i-1} (\mathcal{X}_{s}, \qit) \subset N^{1}H^{2i-1} (\mathcal{X}_{s}, \qit) \]
\item\label{local} The $G_{\cit(S)}$-module $P_{*}H^{2i-1}_{\text{\'et}} (\X_{\overline{\cit(S)}}, \qit_{\ell})$ is irreducible for some (and hence all) $\ell$.
\end{enumerate}
\end{Def}
\begin{Rem} The prototypical case of this notion is in codimension $2$ for families of principally polarized Abelian varieties with sufficiently large monodromy and $P$ an idempotent giving the primitive cohomology in degree $3$. A fully fleshed-out example involving families of products of elliptic curves (but not quite corresponding to primitive cohomology) will be given in the sequel. 
\end{Rem}

\noindent Our interest in primitive correspondences is the following lemma, whose proof is abstracted from \cite{T} pp. 4-5. 
\begin{Lem}\label{kill-tors} Let $f: \X \to S$ be a family as in Notation \ref{notation} and fix a prime $\ell$. Suppose that for the very general $s \in S(\cit)$, we have $N^{1}H^{3} (\mathcal{X}_{s}, \zit/\ell^{r}) \neq H^{3} (\mathcal{X}_{s}, \zit/\ell^{r})$ for $r>0$.
Finally, let $P \in CH^{d} (\X \times_{S} \X)$ be a primitive correspondence in codimension $2$. Then, there exists some integer $M \geq 0$ for which 
\[ \ell^{M}\cdot P_{*}CH^{2} (\mathcal{X}_{s})[\ell^{\infty}] = 0\]
for the very general $s \in S(\cit)$.
\begin{proof} Since $N^{1}H^{3} (\mathcal{X}_{s}, \zit/\ell^{r}) \neq H^{3} (\mathcal{X}_{s}, \zit/\ell^{r})$ for $r >0$, taking the limit over $r$ and tensoring with $\qit_{\ell}$ shows that  the inclusion
\[ B:= \mathop{\lim_{\longleftarrow}}_{r >0} N^{1}H^{3} (\mathcal{X}_{s}, \zit/\ell^{r}) \otimes_{\zit_{\ell}} \qit_{\ell} \to H^{3} (\mathcal{X}_{s}, \qit) \otimes_{\qit} \qit_{\ell} \]
is not surjective. Moreover, we have $Q_{*}H^{3} (\mathcal{X}_{s}, \qit) \subset N^{1}H^{3} (\mathcal{X}_{s}, \qit)$ and that 
\[P_{*}H^{3} (\mathcal{X}_{s}, \qit) \oplus Q_{*}H^{3} (\mathcal{X}_{s}, \qit) = H^{3} (\mathcal{X}_{s}, \qit)\] 
by assumption. Then, using the fact that the $G_{\cit(S)}$-module 
\[ P_{*}H^{3} (\X_{s}, \qit) \otimes \qit_{\ell} \cong P_{*}H^{3}_{\text{\'et}} (\X_{\overline{\cit(S)}}, \qit_{\ell}) \] is irreducible (also by assumption), it follows as in p. 4 paragraph 3 of loc. cit. that $B= Q_{*}H^{3} (\mathcal{X}_{s}, \qit_{\ell}).$ Arguing as in the last two paragraphs of loc. cit., we then deduce that there is some $M \geq 0$
\begin{equation} \ell^{M}\cdot P_{*}(N^{1}H^{3} (\mathcal{X}_{s}, \zit/\ell^{r})) = 0 \label{goal-kill}\end{equation} 
for all $r > 0$. Now, consider the Bloch cycle class map \cite{B}
\[ CH^{2} (\mathcal{X}_{s})[\ell^{\infty}] \to H^{3} (\mathcal{X}_{s}, \qit_{\ell}/\zit_{\ell}(2))\]
By the Merkurjev-Suslin theorem \cite{MS}, this map is injective and its image is 
\[N^{1}H^{3} (\mathcal{X}_{s}, \qit_{\ell}/\zit_{\ell}(2)).\] 
Combining this with (\ref{goal-kill}) gives the desired result.
\end{proof}
\end{Lem}
\begin{Rem} As noted in \cite{T}, we make the observation that (in general) it is possible to have $N^{1}H^{3} (\mathcal{X}_{s}, \qit_{\ell}) \neq B$ (see, for instance, \cite{S4} after Theorem 0.4). In particular, knowing that $N^{1}H^{3} (\mathcal{X}_{s}, \qit_{\ell}) \neq H^{3} (\mathcal{X}_{s}, \qit_{\ell})$ is not sufficient in the above proof, highlighting the need for the Bloch-Esnault result. 
\end{Rem}
\noindent As an application, fix a prime $\ell$ and let $P$ be a primitive correspondence in codimension $2$. Now, let
\[ CH^{2} (\mathcal{X})_{0} := \text{ker } \{ CH^{2} (\mathcal{X}) \xrightarrow{c^{2}} H^{4} (\mathcal{X}, \zit(2)) \xrightarrow{d_{2}^{0,4}} H^{0} (S, R^{4}f_{*}\zit(2)) \} \]
where $ c^{2}$ is the cycle class map and $d_{2}^{0,4}$ is the differential in the Leray spectral sequence
\[ E_{2}^{p,q}= H^{p} (S, R^{q}f_{*}\zit(n)) \Rightarrow H^{p+q} (\X,\zit(n)). \] Then, there is an Abel-Jacobi map:
\begin{equation} CH^{2} (\X)_{0} \to H^{1} (S, R^{3}f_{*}\zit(2))\label{AJ} \end{equation}
arising also from this spectral sequence. We also consider the Galois cohomological version of the Abel-Jacobi map. Namely, for any finite extension $E/\cit(S)$ 
\begin{equation} CH^{2}_{hom} (\X_{E}) \to  H^{1} (G_{E}, H^{3}_{\text{\'et}} (\X_{\overline{\cit(S)}}, \zit/\ell^{r}(2)))\label{abel}\end{equation}
arising from the Hochschild-Serre spectral sequence 
\[ E_{2}^{p,q}= H^{p} (G_{E}, H^{q}_{\text{\'et}} (Y_{\overline{E}}, \zit/\ell^{r}(n))) \Rightarrow H^{p+q}_{\text{\'et}} (Y_{E},\zit/\ell^{r}(n)) \]
for $Y$ a smooth projective variety over $E$. By taking an inverse limit of (\ref{abel}) over $r$, we obtain
\begin{equation} CH^{2}_{hom} (\X_{E}) \to  \mathop{\lim_{\longleftarrow}}_{r \geq 1} H^{1} (G_{E}, H^{3}_{\text{\'et}} (\X_{\overline{\cit(S)}}, \zit/\ell^{r}(2))) \label{abel-2}\end{equation}
\begin{Def} With the above notation, we say that $\gamma \in CH^{2} (\mathcal{X})_{0}$ is {\em primitive with respect to $P$} if $P_{*}(\gamma)$ is non-torsion under (\ref{AJ}).\end{Def}
\noindent There is then the following straightforward observation:
\begin{Lem}\label{AJ-lem} If $\gamma$ is primitive, the image of $P_{*}(\gamma)$ under (\ref{abel-2}) is non-torsion.
\begin{proof} By the Gysin sequence, it follows that the restriction map
\[H^{1} (S, R^{3}f_{*}\zit/\ell^{r}(2))) \to  H^{1} (U, R^{3}f_{U*}\zit/\ell^{r}(2))) \]
is injective for any Zariski open subset $U \subset S$ and $r\geq 1$. Taking the direct limit over all $U$ (and using the comparison isomorphism between \'etale cohomology and singular cohomology) and the inverse limit over $r$, it follows that the image of $P_{*}(\gamma)$ under (\ref{abel-2}) is non-torsion for $E = \cit(S)$. The case of general $E$ then follows from the fact that the kernel of the base extension map:
\[\mathop{\lim_{\longleftarrow}}_{r \geq 1} H^{1} (G_{\cit(S)}, H^{3}_{\text{\'et}} (\X_{\overline{\cit(S)}}, \zit/\ell^{r}(2))) \to \mathop{\lim_{\longleftarrow}}_{r \geq 1} H^{1} (G_{E}, H^{3}_{\text{\'et}} (\X_{\overline{\cit(S)}}, \zit/\ell^{r}(2)))\]
is torsion by a standard transfer argument.

\end{proof}
\end{Lem}

\section{Nori's recipe for infinitude}
\subsection*{Fields of definition}
\noindent Nori's isogeny technique \cite{N} (which was adapted by Totaro to the mod $\ell$ case) can be used quite generally to produce infinitely many distinct cycles in the situation described here. As we will show, this method works so long as ``field of definition" of the cycles (defined below) being used is a ramified extension of the function field of a modular variety. By using isogenies, one can generate infinitely many distinct cycles since the isogenies end up permuting the branch locus.
\begin{Def} Let $L$ be a field (of characteristic $0$, for simplicity) and $n, j$ positive integers, $X$ a smooth projective variety over $L$ and 
\[ \gamma \in CH^{j} (X_{\overline{L}}) \otimes \zit/n \]
The {\em (minimum) field of definition of $\gamma$ over $L$} is the fixed field $K$ of $\overline{L}$ by the stabilizer of $\gamma \in  CH^{j} (X_{\overline{L}}) \otimes \zit/n$ in $G_{L}$. When there is no confusion, we will omit ``over $L$."
\end{Def}
\noindent We observe that stabilizer of $\gamma \in  CH^{j} (X_{\overline{L}}) \otimes \zit/n$ is a closed subgroup of finite index in $G_{L}$. Indeed, $\gamma$ lies in the image of the map
\begin{equation} CH^{j} (X_{F}) \otimes \zit/n \to CH^{j} (X_{\overline{L}}) \otimes \zit/n \label{extens} \end{equation}
for some finite extension $F/L$. It follows that the stabilizer of $\gamma$ contains $G_{F}$ and hence that the stabilizer is a subgroup of $G_{L}$ of finite index. In particular, this means that the field extension $K/L$ is finite.
\begin{Rem}\label{note} The choice of nomenclature here is completely ad hoc and may be a bit deceptive. We note in particular that the field of definition need not be the minimum field $F$ for which $\gamma$ lies in the image (\ref{extens}); if $\gamma$ lies in the image of (\ref{extens}), then $F$ a priori only contains the field of definition of $\gamma$.
\end{Rem}
\begin{Set}\label{setup} Now, suppose that $A_{1}, \ldots A_{m}$ are very general complex Abelian varieties of dimensions $g_{1}, \ldots g_{m}$ for $m\geq 2$ and set $g:= g_{1}+ \ldots +g_{m}$. Suppose that $f_{k}: \mathcal{X}_{k} \to S$ are families as in Notation \ref{notation} for which the geometric generic fibers are isomorphic to $A_{k}$. Then, the geometric generic fiber of the family
\[ f: \X:= \X_{1} \times_{S} \ldots \times_{S} \X_{m} \to S \]
is isomorphic as a scheme over $\cit \cong \overline{\cit(S)}$ to the product
\[ A:= A_{1} \times \ldots \times A_{m} \]
Additionally, we will denote by
\[\A(N) := \A_{g_{1}}(N) \times \ldots \times \A_{g_{m}}(N) \]
where $\A_{g}(N)$ denotes the fine moduli scheme of principally polarized Abelian varieties of dimension $g$ with full level-$N$ structure for $N\geq 3$. Finally, let $\X_{g_{k}} (N) \to \A_{g_{k}}(N)$ denote the corresponding universal Abelian variety and set
\[ \X(N):= \X_{g_{1}} (N) \times \ldots \times \X_{g_{m}} (N) \]
We will assume for simplicity that $S$ has the same dimension as the product of coarse moduli spaces:
\[\A:= \A_{g_{1}} \times \ldots \times \A_{g_{m}} \]
By the assumption on the geometric generic fiber of $f$, it follows that $\overline{\cit(S)} = \overline{\cit(\A)}$. 
\end{Set}
\noindent Let $\gamma \in CH^{2} (\X)$ be a primitive cycle with respect to some primitive correspondence $P \in CH^{g} (\X \times_{S} \X)$, as in the previous section. Moreover, as a matter of convenience for the proof of Proposition \ref{longer}, we impose the following condition:
\begin{Ass}\label{ass} There exists $P_{N} \in CH^{g} (\X(N) \times_{\A(N)} \X(N))$ that coincides with $P$ in the group
\[ CH^{g} (\X(N)_{\overline{\cit(\A(N))}} \times \X(N)_{\overline{\cit(\A(N))}}) = CH^{g} (\X_{\overline{\cit(S)}} \times \X_{\overline{\cit(S)}}) \]
and such that for $s \in S(\cit)$, there are some $j_{1}, \ldots j_{m}$, at least two of which are non-zero, for which
\begin{equation} P_{N*}H^{3} (\X_{s}, \zit) \subset H^{j_{1}} (\X_{1, s}, \zit)\otimes \ldots \otimes  H^{j_{m}} (\X_{m, s}, \zit)\label{two}\end{equation}
\end{Ass}
\begin{Rem} In the above assumption, (\ref{two}) follows from \ref{local} in the definition of primitive correspondence and the K\"unneth theorem. So the condition being imposed in Assumption \ref{ass} is on the number of non-zero $j_{k}$ degrees.
\end{Rem}
\noindent Further, let $\ell$ be a prime, fix $N \geq 3$ and view
\[ \gamma \in CH^{2} (\X_{\overline{\cit(\A(N))}}) \otimes \zit/\ell^{r}, \]
and we allow $r$ to vary over the positive integers. We let $K_{r}$ denote the field of definition of $\gamma$ over $\cit(\A(N))$. We note that a priori that $K_{r}$ could depend on $r$, so we will need the lemma below. (It also can and does depend on $\ell$ and $N$, but since these are fixed, we suppress these from the notation for now.)
\begin{Lem}\label{next} For $r$ sufficiently large, the field of definition $K_{r}$ of $\gamma$ over $\cit(\A(N))$ does not depend on $r$.
\begin{proof}[Proof of Lemma] Note that for all $r \geq 0$, we have $K_{r} \subset E$ for any finite extension $E/\cit(\A(N))$ for which $\gamma$ lies in the image of the natural map
 \begin{equation} CH^{2} (\X_{E}) \to CH^{2} (\X_{\overline{\cit(\A(N))}})\label{natural-map} \end{equation} 
Fix such an extension; then, it is an easy exercise to show that $K_{r} \subset K_{r'}$ whenever $r \leq r'$. Since $E$ is a finite extension of $\cit(\A(N))$, it follows that the ascending chain of $K_{r}$'s must eventually stabilize. Hence, the result.
\end{proof}
\end{Lem}
\noindent Henceforth, we will assume that $r$ is sufficiently large so that the field of definition of $\gamma$ will not depend on $r$. Now, it is possible for the field of definition of $\gamma$ to be what we call {\em a congruence subfield over $\A(N)$} (or simply {\em a congruence subfield}), which we define to be a function field of a finite \'etale cover of some $\A(N)$. The results below show that this is not the case.

\begin{Prop}\label{longer} Let $E/\cit(\A(N))$ be a Galois extension for which $\gamma \in CH^{2} (\X_{E})$ and that $F$ is a congruence subfield over $\A(N)$ with $\cit(\A(N)) \subset F \subset E$. Then, there is some $g \in Gal(E/F)$ for which:
\[[g(P_{N*}(\gamma))] \neq [P_{N*}(\gamma)] \in H^{1}(G_{E}, H^{3}_{\text{\'et}} (\X_{\overline{\cit(\A(N))}}, \zit/\ell^{r}(2))) \]
for some $r \gg 0$, where $[\cdot]$ denotes the Abel-Jacobi image of a cycle.
\begin{proof} Suppose by way of contradiction that there is some congruence subfield $F \subset E$ over $\A(N)$ such that for each $g \in Gal(E/F)$:
\begin{equation}
\label{inv-2}[g(P_{N*}(\gamma))] = [P_{N*}(\gamma)] \in H^{1}(G_{E}, H^{3}_{\text{\'et}} (\X_{\overline{\cit(\A(N))}}, \zit/\ell^{r}(2)))
\end{equation}
for $r$ sufficiently large. Since $Gal(E/F)$ is finite, it follows that there is some $r_{0}$ such that (\ref{inv-2}) holds for each $g \in Gal(E/F)$ and $r \geq r_{0}$. Then, using the results of \cite{J}, it follows that
\[ H^{1} (G_{E}, H^{3}_{\text{\'et}} (\X_{\overline{\cit(\A(N))}}, \zit_{\ell}(2))) \cong \mathop{\lim_{\longleftarrow}}_{r \geq r_{0}} H^{1} (G_{E}, H^{3}_{\text{\'et}} (\X_{\overline{\cit(\A(N))}}, \zit/\ell^{r}(2))) \]
where the left-hand side denotes continuous Galois cohomology group defined in \cite{J}. It follows that
\begin{equation} [P_{N*}(\gamma)] \in H^{1} (G_{E}, P_{N*}H^{3}_{\text{\'et}} (\X_{\overline{\cit(\A(N))}}, \zit_{\ell}(2)))^{Gal(E/F)} \label{desc} \end{equation}
Using Assumption \ref{ass}, we have
\begin{equation} P_{N*}H^{3}_{\text{\'et}} (\X_{\overline{\cit(\A(N))}}, \qit_{\ell}(2)) = P_{*}H^{3}_{\text{\'et}} (\X_{\overline{\cit(S)}}, \qit_{\ell}(2))\label{prim-eq} \end{equation}
which is irreducible as a $G_{E}$-module by the definition of primitive correspondence. (By definition, (\ref{prim-eq}) is irreducible as a $G_{\cit(S)}$-module, so it is also irreducible as a module over any finite index subgroup of $G_{\cit(S)}$.) Thus, it follows that 
\begin{equation} P_{N*}H^{3}_{\text{\'et}} (\X_{\overline{\cit(\A(N))}}, \zit_{\ell}(2))^{G_{E}} = 0\label{vanishes} \end{equation}
Then, there is the Hochschild-Serre spectral sequence (see \cite{J} p. 208), which gives a short exact sequence:
\[ \begin{split}H^{1} (G', P_{N*}H^{3}_{\text{\'et}} (\X_{\overline{\cit(\A(N)}}, \zit_{\ell}(2))^{G_{E}}) & \hookrightarrow H^{1} (G_{F}, P_{N*}H^{3}_{\text{\'et}} (\X_{\overline{\cit(\A(N))}}, \zit_{\ell}(2)))\\ & \to H^{1} (G_{E}, P_{N*}H^{3}_{\text{\'et}} (\X_{\overline{\cit(\A(N))}}, \zit_{\ell}(2)))^{G'} \\ & \to H^{2} (G', P_{N*}H^{3}_{\text{\'et}} (\X_{\overline{\cit(\A(N))}}, \zit_{\ell}(2))^{G_{E}}) \end{split}\]
where $G'= Gal(E/F)$. Because of (\ref{vanishes}) we deduce that
\begin{equation}\label{equality}
H^{1} (G_{F}, P_{N*}H^{3}_{\text{\'et}} (\X_{\overline{\cit(\A(N))}}, \zit_{\ell}(2))) = H^{1} (G_{E}, P_{N*}H^{3}_{\text{\'et}} (\X_{\overline{\cit(\A(N))}}, \zit_{\ell}(2)))^{G'}
\end{equation}
(We note that this argument also appears in \cite{T} Lemma 2.3). By Lemma \ref{AJ-lem}, we obtain a non-torsion cycle in
\[ H^{1} (G_{F}, P_{N*}H^{3}_{\text{\'et}} (\X_{\overline{\cit(\A(N))}}, \zit_{\ell}(2))) \]
In fact, we can spread this non-torsion cycle out in the following way. Let $\X_{g_{i}} (N) \to \A_{g_{i}} (N)$ be the universal family over $\A_{g_{i}} (N)$ and set 
\[ \X(N) = \X_{g_{1}} (N) \times \ldots \times  \X_{g_{m}} (N).\]
Moreover, let $\Y$ be the normalization of $\A(N)$ in $F$. Since $F$ is a congruence subfield, $\Y$ is a finite \'etale cover of $\A(N)$. By slight abuse of notation, we will denote $f_{\Y}: \X_{\Y}:= X(N) \times_{\A(N)} \Y \to \Y$. Then, we have the following claim:
\begin{claim} $[P_{N*}(\gamma)]$ lies in the image of the natural map
\begin{equation} H^{1}_{\text{\'et}} (\Y, R^{3}f_{\Y*} \zit/\ell^{r}(2)) \to H^{1}(G_{F}, H^{3}_{\text{\'et}} (\X_{\overline{\cit(\A(N))}}, \zit/\ell^{r}(2)))\label{restriction} \end{equation}
\begin{proof}[Proof of claim] Note that
\[ H^{1}(G_{F}, H^{3}_{\text{\'et}} (\X_{\overline{\cit(\A(N))}}, \zit/\ell^{r}(2))) = \mathop{\lim_{\longrightarrow}}_{U} H^{1}_{\text{\'et}} (U, R^{3}f_{U*} \zit/\ell^{r}(2))\]
where $U$ ranges over all open subsets of $\Y$ whose complement is a divisor. This means that there exists some $U$ for which $[\gamma']:=[P_{N*}(\gamma)]$ lies in the image of 
\[ H^{1}_{\text{\'et}} (U, R^{3}f_{U*} \zit/\ell^{r}(2)) \to H^{1}(G_{F}, H^{3}_{\text{\'et}} (\X_{\overline{\cit(\A(N))}}, \zit/\ell^{r}(2)))\]
To show that $[\gamma']$ lies in the image of (\ref{restriction}), we use the Gysin exact sequence to reduce the problem to showing that $[\gamma']$ lies in
\begin{equation} \bigcap_{D} \text{ker} \{H^{1}_{\text{\'et}} (U, R^{3}f_{U} \zit/\ell^{r}(2)) \xrightarrow{res_{D}} H^{0} (G_{\cit(D)}, H^{3}_{\text{\'et}} (\X_{\overline{\cit(D)}}, \zit/\ell^{r}(1))) \}\label{residue}
\end{equation}
where $D$ ranges over all the irreducible components of the complement of $U$ in $\Y$ and $res_{D}$ denotes the residue map along $D$. Now, let $\Y'$ denote the normalization of $\Y$ in $E$ and denote the corresponding finite, flat morphism $\pi: \Y' \to \Y$. We can spread out $\gamma'$ to obtain a cycle in $CH^{2} (\X_{\Y} \times_{\Y} \Y') \otimes \zit/\ell^{r}$; necessarily, the image of this cycle under (\ref{abel}) already lies in the kernel of all the residue maps
\[H^{1}(G_{E}, H^{3}_{\text{\'et}} (\X_{\overline{\cit(\A(N))}}, \zit/\ell^{r}(2))) \xrightarrow{res_{D'}} H^{0} (G_{\cit(D')}, H^{3}_{\text{\'et}} (\X_{\overline{\cit(D')}}, \zit/\ell^{r}(1)))
\] 
where $D'$ ranges over the irreducible components of the complement of $\pi^{-1}(U)$ in $\Y'$. Since $\pi$ is finite and flat, $\pi^{-1}(D)$ is a non-empty divisor for each $D$ in $\Y \setminus U$. Letting $D'$ be an irreducible component of $\pi^{-1}(D)$; then, there is the obvious commutative diagram
\[ \begin{tikzcd}
H^{1}(G_{F}, H^{3}_{\text{\'et}} (\X_{\overline{\cit(\A(N))}}, \zit/\ell^{r}(2))) \arrow{r}{res_{D}} \arrow{d}{\pi^{*}} & H^{0} (G_{\cit(D)}, H^{3}_{\text{\'et}} (\X_{\overline{\cit(D)}}, \zit/\ell^{r}(1)))\arrow[hook]{d}{\pi^{*}}\\
H^{1}(G_{E}, H^{3}_{\text{\'et}} (\X_{\overline{\cit(\A(N))}}, \zit/\ell^{r}(2))) \arrow{r}{\oplus res_{D'}} & \displaystyle \bigoplus_{D'}H^{0} (G_{\cit(D')}, H^{3}_{\text{\'et}} (\X_{\overline{\cit(D')}}, \zit/\ell^{r}(1)))
\end{tikzcd}\]
where the direct sum in the bottom row is over all irreducible components of $\pi^{-1} (D)$. It follows that $[\gamma']$ lies in (\ref{residue}), as desired.
\end{proof}
\end{claim}
\noindent By taking an inverse limit over $r$, it now follows that there is a non-torsion cycle in the group
\[ \mathop{\lim_{\longleftarrow}}_{r \geq r_{0}} H^{1}_{\text{\'et}} (\Y, P_{N*}R^{3}f_{\Y*} \zit/\ell^{r}(2)) \]
We will obtain a contradiction to this by showing that the singular cohomology group
\begin{equation}\label{mid-goal} H^{1} (\Y, P_{N*}R^{3}f_{\Y*} \qit(2)) = 0
\end{equation}
\noindent Towards (\ref{mid-goal}), we first have the following lemma:
\begin{Lem}\label{neat} There is some finite \'etale cover of $\Y$ of the form $\Y_{1} \times \ldots \times \Y_{m}$, where each $\Y_{i}$ is a finite \'etale cover of $\A_{g_{i}}(N)$.
\begin{proof}[Proof of Lemma] Observe that
\[ \pi_{1}(\Y) \leq \pi_{1}(\A(N)) = \prod_{i=1}^{m} \Gamma_{g_{i}} (N)\]
where $\Gamma_{g_{i}} (N) \leq Sp_{2g_{i}} (\zit)$ is the full level-$N$ congruence subgroup of the symplectic group. It suffices to show that there exist subgroups of finite index
\[ \Gamma_{i} \leq Sp_{2g_{i}} (\zit) \]
such that $\Gamma := \prod_{i=1}^{m} \Gamma_{i} \leq \pi_{1}(Y)$ is a subgroup of finite index.
To see this, note that $\pi_{1}(Y)$ is a discrete subgroup of the semi-simple Lie group 
\[\prod_{i=1}^{m} Sp_{2g_{i}} (\rit)  \]
Since none of the factors of this group are compact, it follows by \cite{WM} Proposition 4.3.3 (a consequence of the Borel Density theorem) that the desired $\Gamma_{i}$ exist. 
\end{proof}
\end{Lem}
\noindent We obtain (\ref{mid-goal}) with the following result, completing the proof of Proposition \ref{long}.
\begin{Lem} $H^{1} (\Y, P_{N*}R^{3}f_{\Y*}\qit) = 0$.
\begin{proof}[Proof of Claim] By Lemma \ref{neat}, we can assume without loss of generality that 
\[\Y = \Y_{1} \times \ldots \times \Y_{m}\] 
for $\Y_{i}$ a finite \'etale cover of $\A_{g_{i}}(N)$. Then, we have that
\[ \X_{\Y} = \X(N) \times_{\A(N)} \Y \cong \prod_{k=1}^{m} \X_{k, \Y}(N) \]
as schemes over $\Y$, where $\X_{k, \Y}(N) = \X_{g_{k}}(N) \times_{\A_{g_{k}} (N)} \Y_{k}$. Now, let $\pi_{k}:\Y \to \Y_{k}$ and $f_{k, \Y}: \X_{k, \Y}(N) \to\Y_{k}$ be the projections. By Assumption \ref{ass}, we then have
\[P_{N*}R^{3}f_{\Y*}\qit \subset  \pi_{1}^{*}R^{j_{1}}f_{1,\Y*}\qit \otimes \ldots \otimes \pi_{m}^{*}R^{j_{m}}f_{m,\Y*}\qit\]
for $j_{1}+ \ldots + j_{m} =3$ with at least two non-zero $j_{k}$'s. By the K\"unneth theorem, this means that
\begin{equation} H^{1} (\Y, P_{N*}R^{3}f_{\Y*}\qit)\label{equation} \end{equation}
is a subspace of a direct sum of vector spaces of the form:
\begin{equation} H^{i_{1}} (\Y_{1}, R^{j_{1}}f_{1,\Y*}\qit) \otimes \ldots \otimes H^{i_{m}} (\Y_{m}, R^{j_{m}}f_{m,\Y*}\qit)\label{summands} \end{equation}
where exactly one of the indices $i_{1}, \ldots, i_{m}$ is $1$ and the remainder are $0$. Since at least two of the $j_{k}$ indices are non-zero, there are only the following two possibilities (up to permuting the $j_{k}$'s):
\begin{equation} (j_{1}, \ldots , j_{m}) = (2, 1, 0, \ldots 0), \ (1, 1, 1, 0, \ldots 0) \label{possib} \end{equation}
Now, since $\Y_{k}$ is a finite \'etale cover of $\A_{g_{k}}(N)$, it follows that 
\begin{equation} H^{0} (\Y_{k}, R^{1}f_{k,\Y*}\qit) = 0\label{above}\end{equation}
For this same reason, we have that $H^{0} (\Y_{k}, R^{2}f_{k,\Y*}\qit)$ is generated by the class of the corresponding polarization $\Theta_{k}$ (ignoring twists). It follows that any summand in (\ref{summands}) for which no $j_{k}$ is greater than $1$ automatically vanishes. Similarly, any summand in (\ref{summands}) possessing the tensor factor $H^{1} (\Y_{k}, R^{2}f_{k,\Y*}\qit)$ must also possess a tensor factor of the form $H^{0} (\Y_{k}, R^{1}f_{k, \Y*}\qit)$ and must therefore vanish. Finally, any summand in (\ref{summands}) possessing the tensor factor $H^{0} (\Y_{k}, R^{2}f_{k,\Y*}\qit)$ must also possess a tensor factor of the form $H^{i_{k}} (\Y_{k}, R^{1}f_{k, \Y*}\qit)$. If $i_{k} = 0$, this summand again vanishes. If $i_{k}=1$, then the summand is given by
\begin{equation} H^{1} (\Y_{k}, R^{1}f_{k, \Y*}\qit) \otimes H^{0} (\Y_{k'}, R^{2}f_{k',\Y*}\qit) \subset  H^{1} (\Y, \pi_{k'}^{*}\Theta_{k'}\cdot R^{1}f_{\Y*}\qit)\end{equation}
Since $\pi_{k'}^{*}\Theta_{k'}\cdot H^{1} (\X_{s}, \qit) \subset N^{1}H^{1} (\X_{s}, \qit)$, $P_{N*} = P_{*}$ annihilates this summand by \S3. Thus, $P_{N*}$ annihilates all possible summands in (\ref{summands}), and so the desired result follows.
 \end{proof}
\end{Lem}

\end{proof}
\end{Prop}
\begin{Prop}\label{pre-cor} Let $E/\cit(\A(N))$ be a Galois extension for which $\gamma \in CH^{2} (\X_{E})$ and that $F$ is a congruence subfield over $\A(N)$ with $\cit(\A(N)) \subset F \subset E$. Then, there is some $g \in G_{F}$ for which: 
\begin{equation} g(\gamma)\neq \gamma \in CH^{2} (\X_{\overline{\cit(\A(N))}})/\ell^{r} \label{not-eq}\end{equation}
for $r \gg 0$.
\begin{proof} By Proposition \ref{longer}, there is some $g \in Gal(E/F)$ for which 
\begin{equation} [g(P_{N*}(\gamma))] \neq [P_{N*}(\gamma)] \in H^{1}(G_{E}, H^{3}_{\text{\'et}} (\X_{\overline{\cit(\A(N))}}, \zit/\ell^{r}(2)))\label{AJ-key}\end{equation}
for some $r \gg 0$. Hence, (\ref{AJ-key}) holds for all $r$ sufficiently large; say, for all $r \geq r_{1}$. Now, lift $g \in Gal(E/F)$ to some element $\tilde{g} \in G_{F}$. Then, we set
\[ \gamma_{1} := g(\gamma) - \gamma \in CH^{2} (\X_{E}).  \]
so that the base-extension of $\gamma_{1}$ to $\overline{\cit(\A(N))}$ coincides with
\[ \tilde{g}(\gamma) - \gamma \in CH^{2} (\X_{\overline{\cit(\A(N))}}) \]
(where we abuse notation slightly and identify cycles with their base-extensions). Our goal is to show that for $r \geq r_{1} +M$ (where $M$ is as in Lemma \ref{kill-tors}), we have:
\[ \gamma_{1} \neq 0 \in CH^{2} (\X_{\overline{\cit(\A(N))}})/\ell^{r} \] 
To this end, we suppose by way of contradiction that $\gamma_{1} = \ell^{r}\cdot \gamma_{2}$ for some such $r$ and some $\gamma_{2} \in CH^{2} (\X_{\overline{\cit(\A(N))}})$. Since $\gamma_{1}$ is a cycle on $\X_{E}$, it follows that for every $h \in G_{E}$,
\[ 0 = h(\gamma_{1}) - \gamma_{1} = \ell^{r}\cdot (h(\gamma_{2}) - \gamma_{2}) \in CH^{2} (\X_{\overline{\cit(\A(N))}})\]
Applying $P_{N*}$, it follows from Lemma \ref{kill-tors} that
\[ 0 = \ell^{M}\cdot (h(P_{N*}(\gamma_{2})) - P_{N*}(\gamma_{2})) \] 
for all $h \in G_{E}$. Hence, 
\[\ell^{M}\cdot P_{N*}(\gamma_{2}) \in CH^{2} (\X_{\overline{\cit(\A(N))}})^{G_{E}} \]
Then, for some finite Galois extension $K/E$, we have
\[ \gamma_{3}:=\ell^{M}\cdot P_{N*}(\gamma_{2}) \in CH^{2} (\X_{K})^{Gal(K/E)} \]
The argument at the beginning of the proof of Proposition \ref{longer} then gives 
\[ H^{1}(G_{K}, P_{N*}H^{3}_{\text{\'et}} (\X_{\overline{\cit(\A(N))}}, \zit_{\ell}(2)))^{Gal(K/E)} = H^{1}(G_{E}, P_{N*}H^{3}_{\text{\'et}} (\X_{\overline{\cit(\A(N))}}, \zit_{\ell}(2)))\]
Thus, taking the image under the Abel-Jacobi map, we have
\[ [P_{N*}(\gamma_{1})] = \ell^{r-M}\cdot [\gamma_{3}]  \in H^{1}(G_{E}, P_{N*}H^{3}_{\text{\'et}} (\X_{\overline{\cit(\A(N))}}, \zit_{\ell}(2))) \]
or, equivalently,
\[ [P_{N*}(\gamma_{1})] = 0 \in H^{1}(G_{E}, P_{N*}H^{3}_{\text{\'et}} (\X_{\overline{\cit(\A(N))}}, \zit/\ell^{r-M}(2))) \]
This gives the desired contradiction since we have
\[ [P_{N*}(\gamma_{1})] \neq 0 \in H^{1}(G_{E}, H^{3}_{\text{\'et}} (\X_{\overline{\cit(\A(N))}}, \zit/\ell^{s}(2)))\]
for $s \geq r_{1}$ as a result of (\ref{AJ-key}).
\end{proof}
\end{Prop}
\begin{Cor} \label{long} With the conditions of Set-up \ref{setup} and Assumption \ref{ass}, the field of definition of $\gamma$ over $\cit(\A(N))$ is not a congruence subfield over $\A(N)$ for any $N\geq 3$.
\begin{proof} If one assumes that the field of definition of $\gamma$ is a congruence subfield, then Proposition \ref{pre-cor} gives some $g \in G_{F}$ for which (\ref{not-eq}) holds. This contradicts the assumption that $F$ is the field of definition of $\gamma$.
\end{proof}
\end{Cor}
\noindent Before stating the next corollary, we define the following notion, to which we have already alluded.
\begin{Def} Given an irreducible normal variety $U$ with function field $K$ and $L/K$ a finite extension. Moreover, let $V$ denote the normalization of $U$ in $L$. Then, we say that $L/K$ is {\em a ramified extension over $U$} (or simply {\em a ramified extension}) if the natural (finite) morphism $V \to U$ is ramified.
\end{Def}

\begin{Cor}\label{needed-cor} The field of definition of $\gamma$ over $\cit(\A(N))$ is a ramified extension, $K (N)/\cit(\A(N))$. Moreover, for any positive integer $N'$ divisible by $N$, we have 
\[ K(N') = K(N)\cdot\cit(\A(N')) \] 
the compositum of $K(N)$ and $\cit(\A(N'))$.
\begin{proof} The first statement is just a restatement of Corollary \ref{long}, since $K (N)/\cit(\A(N))$ is a ramified extension over $\A(N)$ $\Leftrightarrow$ $K(N)$ is not a congruence subfield over $\A(N)$ (tautologically). For the second statement, note that by definition $K (N)$ is the field corresponding to the stabilizer in the Galois group $G_{\cit(\A(N))}$ of 
\[ \gamma \in CH^{2} (\X_{\overline{\cit(\A(N))}})\otimes \zit/\ell^{r} \]
On the other hand, the absolute Galois group of the compositum $K(N)\cdot\cit(\A(N'))$ is
\[ G_{\cit(\A(N'))} \cap G_{K(N)} \]
which is precisely the stabilizer of $\gamma$ in $\cit(\A(N'))$. The corresponding field is $K(N')$.
\end{proof}
\end{Cor}
\begin{Rem}\label{definition} When $\gamma$ is the Ceresa cycle on the very general Abelian threefold, the field of definition of $\gamma$ is the function field of $\M_{3}(N)$, the fine moduli scheme of curves of genus $3$ with full level $N$-structure ($N\geq 3$), which is a ramified quadratic extension of the function field of $\cit(\A_{3}(N))$, the fine moduli scheme of Abelian threefolds with full level $N$-structure. (By the Torelli theorem, the extension $\cit(\M_{3}(N))/\cit(\A_{3}(N))$ is ramified along the hyperelliptic locus.) 
\end{Rem}
\subsection*{Nori's method}
\noindent To implement Nori's recipe, we consider the following field
\[L := \mathop{\lim_{\longrightarrow}}_{N} \cit(\A(N))\]
where the limit ranges over all $N$ prime to $\ell$. Also, for each $k=1,\ldots m$, consider the  groups
\[ G_{k} = GSp_{2g_{k}} (\qit), \ M_{k} = M_{2g_{k}} (\zit) \]
and set
\[ G= G_{1} \times \ldots \times G_{m}, \ M = M_{1} \times \ldots \times M_{m}\]
Furthermore, we set
\[ W:= \cit^{2g_{1}} \times \ldots \times \cit^{2g_{m}}, \ \mathfrak{h}:=\mathfrak{h}_{g_{1}} \times \ldots \times \mathfrak{h}_{g_{m}} \]
We observe that $G$ acts naturally on $\mathfrak{h} \times W$, the universal cover of the universal Abelian variety $\mathcal{X}(N) \xrightarrow{f_{N}} \mathcal{A}(N)$, with the center 
\[Z(G) = Z(G_{1}) \times \ldots \times Z(G_{m})\] 
acting trivially. Indeed, each
\[\phi \in G \cap M \]
certainly acts on $\mathfrak{h} \times W$ (see, for instance, \cite{N}), and one can extend this action to all of $G$ by taking the action of each $Z(G_{k})$ (consisting of scalar matrices) to be trivial. This action then descends to an action on universal Abelian varieties; i.e., for each $\phi \in G$ and positive integer $N$ prime to $\ell$, there is some $N' \mid N$ and a commutative diagram:
\[ \begin{CD} \X(N) @>{\phi}>> \X(N')\\
@V{f_{N}}VV @V{f_{N'}}VV\\
\A(N)@>{\phi}>> \A(N')  \end{CD}\]
Taking the inverse limit over $N$ we obtain a commutative diagram:
\begin{equation} \begin{CD} \X_{L} @>{\phi}>> \X_{L}\\
@VVV @VVV\\
\text{Spec } L@>{\phi}>> \text{Spec } L\label{generic}  \end{CD}\end{equation}
Here, we set $\X_{\cit(\A(N))}$ to be the generic fiber of $f_{N}$ and 
\[ \X_{L} := \mathop{\lim_{\longleftarrow}}_{N} \X_{\cit(\A(N))}\] 
where again the limit ranges over $N$ prime to $\ell$. This is an Abelian variety of dimension $g$ over the field $L$, and we observe that by construction, each $\phi \in G$ acts by an isogeny (in the diagram (\ref{generic})). As in Nori's approach, we will need a way of keeping track of all the lifts of a given $\phi \in G$ to an action on the geometric generic fiber:
\[ \mathcal{X}_{\overline{L}} := \mathcal{X}_{L} \times_{L} \overline{L} \]
So, we consider the group 
\[\G := \{ (\psi, \phi)\in  Aut(\text{Spec}(\overline{L})) \times G \ | \restr{\psi}{L} = \restr{\phi}{L} \} \] where the notation $\restr{\phi}{L}$ means the restriction of $\phi \in G$ to its action on $L$. Note that there is a natural short exact sequence
\[ 1 \to G_{L} \to \mathcal{G} \xrightarrow{\pi_{2}} G \to 1 \]
As in op cit. p. 194, one has
\begin{Lem}\label{easy} There is a homomorphism
\[ \rho:\G \to \End(\mathcal{X}_{\overline{L}}) \to \End (CH^{2}(\mathcal{X}_{\overline{L}})\otimes \zit/\ell^{r}) \]
for which the composition 
\[ G_{L} \hookrightarrow \G \xrightarrow{\rho} \End (CH^{2}(\mathcal{X}_{ \overline{L}})\otimes \zit/\ell^{r}) \] is the usual Galois action. 
\end{Lem}
\noindent Let $\End(\mathcal{X}_{\overline{L}})_{\ell}$ be the subgroup of isogenies of degree prime-to-$\ell$, let $\G_{\ell} \leq \G$ be the corresponding subgroup in $\G$ and
\[ G_{\ell} := \pi_{2}(\G_{\ell}) \]
It follows from the argument of \cite{T} Lemma 3.2 that  
\[\rho(\G_{\ell}) \subset \Aut(CH^{2}(\mathcal{X}_{\overline{L}})\otimes \zit/\ell^{r})\]
\indent Now, with the conditions and notation of the previous subsection, consider a primitive cycle $\gamma \in CH^{2} (\X)$ with respect to some primitive correspondence $P$ on $\X$ satisfying Assumption \ref{ass}. Then, let 
\[ G' \leq G_{\ell} \cap (Sp_{2g_{1}} (\qit) \times \ldots \times Sp_{2g_{m}} (\qit)) \]
be a subgroup for which there are
\[\phi_{1}, \phi_{2}, \ldots \in G'\] 
representing infinitely many distinct cosets in $G/\Gamma$ for 
\[ \Gamma := Sp_{2g_{1}} (\zit) \times \ldots \times Sp_{2g_{m}} (\zit). \]
Lift these to $\tilde{\phi}_{1}, \tilde{\phi}_{2}, \ldots \in \G$ via $\mathcal{G}' \xrightarrow{\pi_{2}} G'$. Note by assumption that each $\rho(\tilde{\phi}_{j})$ is an automorphism on
\[CH^{2}(\mathcal{X}_{\overline{L}})\otimes \zit/\ell^{r}\] 
Thus, we have the following:
\begin{Prop}\label{Nori} With the above notation, there exist subgroups $G'$ as above and $\phi_{1}, \phi_{2}, \ldots \in G'$ for which the corresponding set
\[ \{ \rho(\tilde{\phi}_{1})(\gamma), \rho(\tilde{\phi}_{2})(\gamma), \ldots \} \subset CH^{2}(\mathcal{X}_{\overline{L}})\otimes \zit/\ell^{r} \]
is infinite when $r$ is sufficiently large.
\begin{proof} It follows by Corollary \ref{needed-cor} that the field of definition of 
\[ \gamma \in CH^{2} (\X_{\overline{\cit(S)}}) \otimes \zit/\ell^{r} = CH^{2} (\X_{\overline{L}}) \otimes \zit/\ell^{r} \]
over $L$ is the field $K_{\gamma} = \displaystyle \mathop{\lim_{\longrightarrow}}_{N} K_{\gamma}(N),$ where the limit ranges over $N$ prime to $\ell$ and where $K_{\gamma}(N)$ is the field of definition of $\gamma$ over $\cit(\A(N))$. For each such $N$, let 
\[B_{\gamma}(N) \subset \mathcal{A}(N)\]
denote the corresponding branch locus of the extension $K_{\gamma}(N)/\cit(\A(N))$ (by which we mean the branch locus of the finite morphism $\A_{\gamma}(N) \to \A(N)$, where $\A_{\gamma}(N)$ is the normalization of $\A(N)$ in $K_{\gamma}(N)$). We note that for $N \mid N'$, $B_{\gamma}(N')$ is the pull-back of $B_{\gamma}(N)$ along $\A(N') \to \A(N)$. This follows from Corollary \ref{needed-cor} and the fact that $\A(N') \to \A(N)$ is finite \'etale. Thus, the pull-back of $B_{\gamma}(N)$ along $\mathfrak{h} \to \A(N)$ is a closed subset of $\mathfrak{h}$ for the analytic topology and does not depend on $N$. We will call this set the {\em absolute branch locus}, $B_{\gamma} \subset \mathfrak{h}$. \\
\indent Each of the cycles $\gamma_{j}:=\rho(\tilde{\phi}_{j})(\gamma)$ has a field of definition (over $L$) $K_{\gamma_{j}}$, and we denote the absolute Galois groups by $H_{j}:= G_{K_{\gamma_{j}}}$. This latter is the stabilizer of 
\[ \gamma_{j} \in CH^{2} (\X_{\overline{L}}) \otimes \zit/\ell^{r} \] in $G_{L}$. In order to show that $\gamma_{j}$ are all distinct, it will suffice to show that the $H_{j}$ are all distinct. Now, by Lemma \ref{easy} (and a basic property of stabilizer subgroups), we have
\[ H_{j} = \tilde{\phi}_{j}G_{K_{\gamma}}\tilde{\phi}_{j}^{-1} \leq \mathcal{G}\]
and, as a consequence, $K_{\gamma_{j}} = \tilde{\phi}_{j}(K_{\gamma})$. In particular, it follows that $K_{\gamma_{j}}$ is a ramified extension of $L$; as in the first paragraph of the proof, one then has a corresponding absolute branch locus $B_{\gamma_{j}}$ and this will satisfy:
\[B_{\gamma_{j}} = \phi_{j}^{-1}(B_{\gamma}) \subset \mathfrak{h}\]
Hence, it will suffice (as in \cite{N} and \cite{T}) to show that all the $B_{\gamma_{j}}$ are distinct. Since $\phi_{j}$ were chosen to have distinct cosets in $G/\Gamma$, it will suffice to prove the lemma below.
\end{proof}
\end{Prop}
\begin{Lem} Consider the group 
\[\Lambda= \{ \phi \in Sp_{2g_{1}} (\rit) \times \ldots \times Sp_{2g_{m}} (\rit)  \ | \  \phi(B_{\gamma}) = B_{\gamma}\} \] 
Then, there exist subgroups $G'$ for which $G'/G' \cap \Gamma$ is infinite and such that $G' \cap \Lambda \leq \Gamma$.
\begin{proof} Note that since $B_{\gamma}$ is obtained by pulling back along $\mathfrak{h} \to \A(N)$, it is stable under the action of a subgroup of finite index of $\Gamma$. Thus, the corresponding Lie algebra $\text{Lie}(\Lambda)$ is stable under the adjoint action $\Gamma$, which is Zariski dense in
\[ K(\rit):= Sp_{2g_{1}} (\rit) \times \ldots \times Sp_{2g_{m}} (\rit) \]
from which it follows that $\text{Lie}(\Lambda)$ is invariant under the adjoint action of $K(\rit)$. Since $K(\rit)$ is a semi-simple Lie group, it follows that $\text{Lie}(\Lambda)$ is a product of the Lie algebras
\[ \mathfrak{s}\mathfrak{p}_{2g_{1}} (\rit), \ldots, \mathfrak{s}\mathfrak{p}_{2g_{m}} (\rit). \]
Since $B_{\gamma} \neq \mathfrak{h}$, we cannot have $\text{Lie}(\Lambda) =\text{Lie}(K(\rit))$, from which we deduce that there is some $j$ for which
\[ \text{Lie}(\Lambda) \cap \mathfrak{s}\mathfrak{p}_{2g_{j}} (\rit) = 0 \]
where we view $\mathfrak{s}\mathfrak{p}_{2g_{j}} (\rit)$ as a summand of $\text{Lie}(K(\rit))$. Since $\Lambda$ is closed in the analytic topology in $K(\rit)$, it follows that 
\[ \Lambda \cap Sp_{2g_{j}} (\rit)\]
is discrete. By the Borel-Density Theorem, we have 
\[\Lambda \cap Sp_{2g_{j}} (\rit) \leq Sp_{2g_{j}} (\zit) \] 
Thus, we can take any subgroup $G' \leq Sp_{2g_{j}} (\qit) \cap G_{\ell}$ for which $G'/G' \cap \Gamma$ is infinite.
\end{proof}
\end{Lem}

\section{Application}
\noindent We would like to apply the method described in the previous sections to obtain a proof of Theorem \ref{main}. For this we will first need families $f_{k}: \X_{k} \to S$ as in Set-up \ref{setup} for which $m=3$ and for which $g_{1}=g_{2}=g_{3}=1$. In particular, this means that $S$ will have dimension $3$. We also need for the very general fiber of the product of the $f_{k}$'s to be a product of very general elliptic curves. While there are certainly many possible choices, we will consider the following extension of a family of elliptic curves considered in \cite{BST} and \cite{S1}.\\
\indent Let $S \subset \pit(H^{0} (\pit_{\cit}^{2}, \mathcal{O}(2)))$ be an open subset consisting of smooth conics in the family defined by:
\begin{equation} x^{2}+y^{2}+z^{2}+axy+byz+cxz = 0 \label{import} \end{equation}
Then, we consider the universal family of conics over $S$:
\begin{equation} \begin{CD}  \mathcal{Q} @>>> & \pit^{2} \\ @VhVV & \\ S &  \end{CD}\label{fiber}\end{equation}
There is a finite morphism $\phi: \pit^{2} \to \pit^{2}$ defined by
\begin{equation} [x, y, z] \mapsto [x^{2}, y^{2}, z^{2}] \label{phi}\end{equation}
whose Galois group is $G:= \zit/2 \times \zit/2$ and is generated by the involutions:
\[ [x, y, z] \xmapsto{\sigma_{1}} [-x, y, z], \ [x, y, z] \xmapsto{\sigma_{2}} [x, -y, z]  \]
and denote by $\sigma_{3} := \sigma_{1}\circ\sigma_{2} = \sigma_{2}\circ\sigma_{1}$. We then consider the Cartesian product:
\begin{equation}\label{cart-diag} \begin{CD}  \mathcal{C} @>q>> &\mathcal{Q}  \\ @VVV & @VVV\\ \pit^{2} @>\phi>> & \pit^{2} \end{CD}\end{equation}
There is a natural morphism $g: \mathcal{C} \xrightarrow{} S$, and after possibly shrinking $S$, we may assume that all the fibers of $g$ are smooth. Explicitly, the total space $\C$ is (some open subset of) the family of degree $4$ curves in $\pit^{2}$ defined by
\begin{equation} x^{4}+y^{4}+z^{4}+ax^{2}y^{2}+by^{2}z^{2}+cx^{2}z^{2} = 0 \label{3-family}\end{equation}
The involutions above then induce an action on $\mathcal{C} \xrightarrow{g} S$, which gives rise to the quotients:
\[ p_{k}: \C \to \X_{k}:= \mathcal{C}/\sigma_{k}, \ f_{k}: \X_{k} \to S, \ k=1,2,3 \]
where $f_{k} \circ p_{k} = g$. It is easy to see that for the general fiber of $g$, each $\sigma_{k}$ fixes precisely $4$ points; using a Riemann-Hurwitz argument, it follows that the general fiber of $f_{k}$ is smooth of genus 1. Upon shrinking $S$ further, we can assume that all the fibers of $f_{k}$ are genus $1$. Moreover, there are induced quotient maps:
\[ q_{k} : \X_{k} \to \mathcal{Q} = \C/G \]
as well as involutions, $(-1)_{\X_{k}}: \X_{k} \to \X_{k},$ which induce the usual action by $-1$ on the fibers of $f_{k}$ (since the quotient is a family of conics). Then, we consider the fiber product
\[ f:= f_{1} \times_{S} f_{2}  \times_{S} f_{3}: \X:= \X_{1} \times_{S} \X_{2} \times_{S} \X_{3} \to S \]
There is then an action by \[\tilde{G} = \zit/2 \times \zit/2 \times \zit/2\] on $f: \X \to S$ generated by involutions
\[ \tilde{\sigma}_{k} : \X \to \X, \ k=1,2,3 \]
which act by $(-1)_{\X_{k}}$ on the $k^{th}$ factor and the identity on the remaining factors. We also use the notation:
\[ \tilde{\sigma}_{jk} := \tilde{\sigma}_{j}\circ\tilde{\sigma}_{k}, \ \tilde{\sigma}_{123}:= \tilde{\sigma}_{1}\circ\tilde{\sigma}_{2}\circ\tilde{\sigma}_{3} \]
Moreover, there is a morphism:
\[ \rho: \C \xrightarrow{\Delta^{(3)}} \C^{3} := \C \times_{S} \C \times_{S} \C \xrightarrow{p} \X \]
where $p:= p_{1} \times_{S} p_{2} \times_{S} p_{3}$. \\
\indent Now, we need to check that the assumption on very general fibers in Set-up \ref{setup} is satisfied. Namely, we would like to show that the very general fiber of $f$ is the product of $3$ very general elliptic curves. To this end, we need to use $j$-invariants. However, since it is not clear if $h: \mathcal{Q} \to S$ possesses any sections, we need to pass to some generically finite cover of $S$. Indeed, for each $k$ we denote by $B_{k} \subset \mathcal{Q}$ the ramification locus of the degree $2$ quotient map $q_{k}$ and let $B:= B_{1} \cup B_{2} \cup B_{3}$. Explicitly, $B$ is the closed subset of $\mathcal{Q}$ defined by $xyz=0$ and is a union of $3$ quadratic multi-sections over $S$. We consider the degree $8$ multiquadratic cover of $S$, over which all $3$ multi-sections split, which we denote by $\tilde{S}$. An explicit computation shows that:
\[ \cit(a,b,c) = \cit(S) \subset \cit(\tilde{S}) =\cit(S)(\sqrt{a^{2}-4},\sqrt{b^{2}-4}, \sqrt{c^{2}-4})  \]
Setting $\displaystyle \chi(t) = \frac{t-2}{t+2}$, we may view $\tilde{S} \subset \mathbb{A}^{3}$ with coordinates $u_{1} = \sqrt{\chi(a)}$, $u_{2} = \sqrt{\chi(b)}$, $u_{3} = \sqrt{\chi(c)}$. Now, the base change
\[ \tilde{h}: \tilde{\mathcal{Q}}:= \mathcal{Q}\times_{S} \tilde{S} \to \tilde{S} \]
becomes a (trivial) $\pit^{1}$-bundle over $\tilde{S}$. Hence, there are associated $j$-invariant maps:
\[ J_{k}: \tilde{S} \to \mathbb{A}^{1} \]
\begin{Lem} There is some $\psi \in \cit(v,w)$ satisfying $\psi(v,w) = \psi(w,v)$ for which 
\begin{equation} J_{k} = \psi(u_{k}, u_{k+2}), \label{equals} \end{equation}
where the subscripts are taken modulo $3$. Moreover, the product of the $j$-invariant maps $J:= J_{1} \times J_{2} \times J_{3}: \tilde{S} \to \mathbb{A}^{3}$ is generically finite.
\begin{proof} In the case at hand, $q_{k}: \X_{k} \to \mathcal{Q}$ realizes $\X_{k}$ as a degree $2$ cover of a rational curve $\mathcal{Q} \to S$ over $S$. The $j$-invariants $J_{k}$ of these are completely determined by the branch locus $B_{k}$ (or rather $\tilde{B}_{k}:= B_{k} \times_{S} \tilde{S}$). For $k=1$, $B_{1} \subset \mathcal{Q}$ is defined by $yz=0$, which we see from (\ref{import}) depends only on $a$ and $c$ (and hence only on $u_{1}$ and $u_{3}$). In fact, writing $s = (a,b,c) \in S$ (in slight abuse of notation), we have
\[ B_{1} \times_{S} s = \{ [-c\pm\sqrt{c^{2}-4} ,0 ,2], [-a\pm\sqrt{a^{2}-4}, 2, 0]  \} \subset \mathcal{Q}_{s} \subset \pit^{2} \]
From this, it follows that $\psi (v,w) = \psi(w,v)$. Moreover, we see that the cyclic permutations on $S$ that permute $a$, $b$ and $c$ lift to cyclic permutations on $\mathcal{Q}$ that also cyclically permutes $x$, $y$ and $z$. These permutations in turn permute $B_{1}$, $B_{2}$ and $B_{3}$, which gives (\ref{equals}). In order to prove the second statement, note that the above symmetry property gives the equality of partials:
\[ \psi_{v} (w,v) = \psi_{w}(v,w) \]
We then consider the Jacobian derivative of $J$ at $u_{1} = u_{2}=u_{3}=v$:
\[ DJ= \begin{pmatrix} 
\psi_{v} (v,v) & 0 & \psi_{v} (v,v)\\
\psi_{v} (v,v) & \psi_{v} (v,v) & 0\\
0 & \psi_{v} (v,v)  & \psi_{v} (v,v)
\end{pmatrix} \]
Computing the determinant gives $2\psi_{v} (v,v)^{3}$. We'd like to show that does not vanish for every $v$. Indeed, if it were to vanish, then $\psi(v,v) = J_{1}(v,v)$ would be constant. However, from \S 4 of \cite{BST}, there is an explicit equation for the elliptic curves in the family $f_{1}^{-1}(\Delta^{(3)}) \to \Delta^{(3)}:= \{ (a, a, a) \in S \}$, which shows that this family is not isotrivial. Hence, the determinant of $DJ$ is non-zero and $J$ is generically finite, as desired.
\end{proof}
\end{Lem}
\noindent We also need to produce a primitive correspondence 
\[ P \in CH^{3} (\X \times_{S} \X) \]
satisfying Assumption \ref{ass}, as well as some $\gamma \in CH^{2} (\X)$ which is primitive with respect to $P$. To define $P$, retain the notation from earlier in the section and set
\[ \Pi_{k} := \Delta_{\X_{k}} - \Gamma_{\tilde{\sigma}_{k}} \in CH^{1} (\X_{k} \times_{S} \X_{k}) \]
and then set
\[ P:= \Pi_{1} \times_{S} \Pi_{2} \times_{S}\Pi_{3} \in CH^{3} (\X_{1}^{\times_{S} 2} \times_{S} \X_{2}^{\times_{S} 2} \times_{S} \X_{3}^{\times_{S} 2}) \cong CH^{3} (\X \times_{S} \X) \] 
The main cycle is then defined as an analogue of the Ceresa cycle:
\begin{equation} \boxed{\gamma := \rho_{*}([\C]) - \tilde{\sigma}_{123}^{*}\rho_{*}([\C]) \in CH^{2} (\X)} \end{equation}
Note that for $\{i,j,k \} = \{1,2,3 \}$ and $j<k$, we have 
\[ \tilde{\sigma}_{jk}^{*}\rho([\C]) = \rho([\C])\in CH^{2} (\X)\] 
since there is a commutative diagram:
\[ \begin{tikzcd} \C \arrow{r}{\sigma_{i}} \arrow{d}{\rho} & \C \arrow{d}{\rho}\\
\X \arrow{r}{\tilde{\sigma}_{jk}} & \X
\end{tikzcd}\]
by construction. It follows that $\gamma$ is invariant under $\tilde{\sigma}_{jk}$ and anti-invariant under $\tilde{\sigma}_{123}$ and hence is anti-invariant with respect to each $\tilde{\sigma}_{k}$. A computation then shows that
\begin{equation} P_{*}(\gamma) = 8\cdot \label{short-cut}\gamma \in CH^{2} (\X)\end{equation}

\subsection*{Proof that $P$ is primitive} We observe that since $\tilde{\sigma}_{k}$ is an involution for $k=1,2,3$, we have
\[ \begin{split}\Pi_{k}^{\circ 2} = (\Delta_{\X_{k}} - \Gamma_{\tilde{\sigma}_{k}})^{\circ 2} &=  \Delta_{\X_{k}} - 2\Gamma_{\tilde{\sigma}_{k}} + \Gamma_{\tilde{\sigma}_{k}^{2}}\\
&= 2(\Delta_{\X_{k}} - \Gamma_{\tilde{\sigma}_{k}})
\end{split} \]
It follows that $P^{\circ 2} = 8\cdot P$ by a standard computation, which verifies the first condition of the definition. To verify condition \ref{coniveau} note that for $s \in S$, we have
Thus, we see that the complementary correspondence
\[ \Phi_{k} = 2\cdot\Delta_{\X_{k}} - \Pi_{k} = \Delta_{\X_{k}} +\Gamma_{\tilde{\sigma}_{k}}= \prescript{t}{}{\Gamma_{q_{k}}}\circ\Gamma_{q_{k}}  \]
from which it follows that
\[\Phi_{k*}H^{*} (\X_{k,s}, \zit) = q_{k}^{*}q_{k*} H^{*} (\mathcal{Q}_{s}, \zit) \subset N^{1}H^{*} (\X_{k,s}, \zit)\]
Now, we compute
\[ \begin{split} Q & = 8\cdot \Delta_{\X} - \Pi_{1} \times_{S} \Pi_{2} \times_{S} \Pi_{3}\\ 
& = \Phi_{1} \times_{S} \Phi_{2} \times_{S} \Phi_{3} +(\Phi_{1} \times_{S} \Phi_{2} \times_{S} \Pi_{3} + \text{perm.}) + (\Phi_{1} \times_{S} \Pi_{2} \times_{S} \Pi_{3} + \text{perm.})\end{split}\]
from which we deduce that 
\[ Q_{*}H^{*} (\X_{s}, \zit) \subset N^{1}H^{*} (\X_{s}, \zit) \]
by properties of the coniveau filtration (i.e., compatibility with cup product). Additionally, we would like to verify condition \ref{local}. Indeed, note that by construction
 \[ \Pi_{k*}H^{1}(\X_{k,s}, \qit_{\ell}) = H^{1} (\X_{k,s}, \qit_{\ell}) \]
Using the K\"unneth theorem again, it then follows that
\[ P_{*}H^{3}(\X_{s}, \qit_{\ell}) = H^{1} (\X_{1,s}, \qit_{\ell}) \otimes H^{1} (\X_{2,s}, \qit_{\ell}) \otimes H^{1} (\X_{3,s}, \qit_{\ell})  \]
For condition \ref{local}, it suffices to show that $V_{\ell}^{\otimes 3}$ is irreducible as a $G_{\cit(\A(N))}$-module, where $\A(N)=\A_{1}(N)^{\times 3}$, $\A_{1}(N)$ is the fine moduli scheme of elliptic curves with full level-$N$ structure (with $N \geq 3$ and coprime to $\ell$), $\X(N)\to \A_{1}(N)$ the universal elliptic curve and $V_{\ell} =H^{1}_{\text{\'et}} (\X(N)_{ \overline{\A_{1}(N)}}, \qit_{\ell})$. However, it is well-known that the Galois image in this case contains the subgroup $GL(V_{\ell})^{\times 3}$, viewed as a subgroup of $GL(V_{\ell}^{\otimes 3})$. Thus, condition \ref{local} is satisfied. Finally, we observe that by taking $P_{N}$ to be defined in the same way as $P$, Assumption \ref{ass} is also satisfied.
\subsection*{Proof that $\gamma$ is primitive} Certainly, 
\[\restr{\gamma}{\X_{s}} \in CH^{2} (\X_{s})\]
is homologically trivial for $s\in S$. We would like to show that the image of $P_{*}(\gamma)$ under (\ref{abel}) in $H^{1} (S, P_{*}R^{3}f_{*}\zit(2))$ is non-torsion or, equivalently that it is non-zero in $H^{1} (S, P_{*}R^{3}f_{*}\qit(2))$. Now, we observe that
\[ P_{*}R^{j}f_{*}\qit(2) = 0\]
for $j\neq 3$ by construction of $P$. From the functoriality of the Leray spectral sequence with respect to correspondences, it follows that
\[ H^{1} (S, P_{*}R^{3}f_{*}\qit(2)) =P_{*}H^{4} (\X, \qit(2)) \]
Hence, it suffices to show that the image of $P_{*}(\gamma) = 8\cdot\gamma$ (see (\ref{short-cut})) under the cycle class map
\begin{equation}CH^{2} (\X) \otimes \qit \to H^{4} (\X, \qit(2))\label{cycle} \end{equation}
is non-zero. To this end, we have the following:
\begin{Lem} Let $\hat{\X}$ be a smooth compactification of $\X$. Then, algebraic equivalence and homological equivalence agree in $CH^{2}(\hat{\X}) \otimes \qit$.
\begin{proof} By \cite{BS} Theorem 1 (ii), it suffices to prove that $\X$ is rationally connected. Since $\Y := \C \times_{S} \C \times_{S} \C$ dominates $\X$, it suffices to prove this for $\Y$. To this end, consider the map
\[ \mathcal{Z}:= \mathcal{Q} \times_{S} \mathcal{Q} \times_{S} \mathcal{Q} \xrightarrow{\pi} Z:=\pit^{2} \times \pit^{2} \times \pit^{2} \] 
obtained by projecting each $\mathcal{Q} \subset S \times \pit^{2}$ to $\pit^{2}$. For a given $t_{1}, t_{2}, t_{3} \in \pit^{2}$, the fiber
\[ \pi^{-1}((t_{1}, t_{2}, t_{3})) \]
consists of the conics in the family (\ref{import}) which pass through $t_{1},t_{2}, t_{3}$. For the general $(t_{1}, t_{2}, t_{3})$, there is a unique such conic, from which it follows that $\pi$ is birational (not necessarily surjective). From (\ref{cart-diag}) it follows that there is a Cartesian diagram:
\[\begin{tikzcd}
\Y \arrow{r} \arrow{d} & \mathcal{Z} \arrow{d}\\
Z \arrow{r}{\phi \times \phi \times \phi} & Z
\end{tikzcd}
\]
where $\phi$ was defined in (\ref{phi}). It follows that $\mathcal{Z} \to Z$ is also birational; hence, so is $\Y \to Z$, which means $\mathcal{Y}$ is rationally connected. 
\end{proof}
\end{Lem}

\noindent Now, suppose by way of contradiction that $\gamma = 0 \in H^{4} (\X, \qit)$. Then, let $\hat{\X}$ be a smooth compactification and let $D = \hat{\X} \setminus \X$ and spread out $\gamma$ to a cycle in $CH^{2} (\hat{\X}) \otimes \qit$, also denoted by $\gamma$ by abuse of notation. By \cite{V} Lemma 1.1, the Voisin standard conjecture holds for codimension $2$ cycles. This means that there exists some $\gamma' \in CH^{2} (\hat{\X}) \otimes \qit$ supported on $D$ for which 
\[ \gamma' = \gamma = 0 \in  H^{4} (\hat{\X}, \qit)\]
Then, by the previous lemma, $\gamma - \gamma' \in CH^{2} (\hat{\X}) \otimes \qit$ is algebraically equivalent to $0$. Moreover, we note that for $s \in S$, the fiber $\X_{s}$ does not intersect $D$, so 
\[ \restr{\gamma'}{\X_{s}} = 0 \in CH^{2} (\X_{s}) \otimes \qit \]
From this, we deduce that
\[ \restr{\gamma}{\X_{s}} \in CH^{2} (\X_{s}) \otimes \qit \]
is algebraically equivalent to $0$ for $s \in S$. (From this, it would also follow that the Ceresa cycle is algebraically equivalent to $0$ in 
\[ CH^{2} (Pic^{0}(\C_{s})) \otimes \qit \]
for $s \in S$.) When $\C_{s}$ is the Fermat quartic, this contradicts \cite{Bl} Theorem 4.1.

\subsection*{Proof of Theorem \ref{main}} We first prove Theorem \ref{main} in the case that $d=2$ and $n=3$. In this case, we let $\X \xrightarrow{f} S$ be as at the beginning of this section, which satisfies the conditions of Set-up \ref{setup}. Moreover, by the previous two sub-sections, $P$ is a primitive correspondence satisfying Assumption \ref{ass}, and $\gamma$ is primitive with respect to $P$, Proposition \ref{Nori} shows that for very general $s \in S$ (and $r$ sufficiently large),
\begin{equation} CH^{2} (\X_{s}) \otimes \zit/\ell^{r}\label{final-Chow} \end{equation}
is infinite. Finally, as in \cite{T}, this implies that $CH^{2} (\X_{s}) \otimes \zit/\ell$ is also infinite.\\
\indent For $n>3$ and $2\leq d\leq n-1$, note that for very general complex elliptic curves $E_{1}, \ldots E_{n}$, one can select degree $2$ maps $f_{k}: E_{k} \to \pit^{1}$. Then, one has
\[ CH^{d} (E_{1} \times E_{2} \times E_{3} \times \underbrace{\pit^{1} \times \ldots \times \pit^{1}}_{n-3 \text{ times}}) \otimes \zit/\ell \]
is infinite for all $d$ in the above range. Moreover, for $\ell \neq 2$, the natural map
\[ CH^{d} (E_{1} \times E_{2} \times E_{3} \times \underbrace{\pit^{1} \times \ldots \times \pit^{1}}_{n-3 \text{ times}}) \otimes \zit/\ell \to CH^{d} (E_{1} \times \ldots \times E_{n}) \otimes \zit/\ell \]
induced by pull-back along 
\[ \text{id}_{E_{1}} \times \text{id}_{E_{2}}  \times \text{id}_{E_{3}} \times f_{4} \times \ldots \times f_{n}: E_{1} \times \ldots \times E_{n} \to E_{1} \times E_{2} \times E_{3} \times \underbrace{\pit^{1} \times \ldots \times \pit^{1}}_{n-3 \text{ times}} \]
is injective (using the projection formula). This gives the desired infinitude for $\ell \neq 2$. For $\ell = 2$, the same argument works except that one selects odd-degree maps $f_{k}: E_{k} \to \pit^{1}$ instead.

\Addresses

\end{document}